\documentclass{amsart}

\usepackage{amssymb} \usepackage{amsthm} \usepackage{amsmath}
\usepackage[all]{xy} \SelectTips{eu}{}
\usepackage[hidelinks]{hyperref}

\newcommand{\numberseries}{\bfseries}   %Fontseries used for numbering

\newlength{\thmtopspace}                %Space above theorem
\newlength{\thmbotspace}                %Space below theorem
\newlength{\thmheadspace}               %Space after theorem label
\newlength{\thmindent}                  %For indenting

\newtheoremstyle{fixed bf head,slanted body}
                {\thmtopspace}{\thmbotspace}{\slshape}
                {\thmindent}{\bfseries}{.}{\thmheadspace}
                {{\numberseries \thmnumber{#2\;}}\thmname{#1}\thmnote{ (#3)}}

\newtheoremstyle{fixed bf head,upright body}
                {\thmtopspace}{\thmbotspace}{\upshape}
                {\thmindent}{\bfseries}{.}{\thmheadspace}
                {{\numberseries \thmnumber{#2\;}}\thmname{#1}\thmnote{ (#3)}}

\newtheoremstyle{numbered paragraph}
                {\thmtopspace}{\thmbotspace}{\upshape}
                {\thmindent}{\upshape}{}{\thmheadspace}
                {{\numberseries \thmnumber{#2.}}}

\theoremstyle{fixed bf head,slanted body}
\newtheorem{thm}{Theorem}[section]          \newtheorem*{thm*}{Theorem}
\newtheorem{prp}[thm]{Proposition}      \newtheorem*{prp*}{Proposition}
\newtheorem{cor}[thm]{Corollary}        \newtheorem*{cor*}{Corollary}
\newtheorem{lem}[thm]{Lemma}            \newtheorem*{lem*}{Lemma}

\theoremstyle{fixed bf head,upright body}
\newtheorem{rmk}[thm]{Remark}           \newtheorem*{rmk*}{Remark}

\theoremstyle{numbered paragraph}
\newtheorem{ipg}[thm]{}

\newlength{\thmlistleft}        %leftmargin
\newlength{\thmlistright}       %rightmargin
\newlength{\thmlistpartopsep}   %partopsep
\newlength{\thmlisttopsep}      %topsep
\newlength{\thmlistparsep}      %parsep
\newlength{\thmlistitemsep}     %itemsep

\newcounter{eqc} 
\newenvironment{eqc}{\begin{list}{\upshape (\textit{\roman{eqc}})}%
    {\usecounter{eqc}%
      \setlength{\leftmargin}{\thmlistleft}%
      \setlength{\labelwidth}{\thmlistleft}%
      \setlength{\rightmargin}{\thmlistright}%
      \setlength{\partopsep}{\thmlistpartopsep}%
      \setlength{\topsep}{\thmlisttopsep}%
      \setlength{\parsep}{\thmlistparsep}%
      \setlength{\itemsep}{\thmlistitemsep}}}%
  {\end{list}}%

\newcommand{\eqclbl}[1]{{\upshape(\textit{#1})}}

\newcounter{prt}
\newenvironment{prt}{\begin{list}{\upshape (\alph{prt})}%
    {\usecounter{prt}%
      \setlength{\leftmargin}{\thmlistleft}%
      \setlength{\labelwidth}{\thmlistleft}%
      \setlength{\rightmargin}{\thmlistright}%
      \setlength{\partopsep}{\thmlistpartopsep}%
      \setlength{\topsep}{\thmlisttopsep}%
      \setlength{\parsep}{\thmlistparsep}%
      \setlength{\itemsep}{\thmlistitemsep}}}%
  {\end{list}}%

\newcounter{rqm}
\newenvironment{rqm}{\begin{list}{\upshape (\arabic{rqm})}%
    {\usecounter{rqm}%
      \setlength{\leftmargin}{\thmlistleft}%
      \setlength{\labelwidth}{\thmlistleft}%
      \setlength{\rightmargin}{\thmlistright}%
      \setlength{\partopsep}{\thmlistpartopsep}%
      \setlength{\topsep}{\thmlisttopsep}%
      \setlength{\parsep}{\thmlistparsep}%
      \setlength{\itemsep}{\thmlistitemsep}}}%
  {\end{list}}%

\newenvironment{itemlist}{\nopagebreak \begin{list}{$\bullet$}%
    {\setlength{\leftmargin}{\thmlistleft}%
      \setlength{\labelwidth}{\thmlistleft}%
      \setlength{\rightmargin}{\thmlistright}%
      \setlength{\partopsep}{\thmlistpartopsep}%
      \setlength{\topsep}{\thmlisttopsep}%
      \setlength{\parsep}{\thmlistparsep}%
      \setlength{\itemsep}{\thmlistitemsep}}}%
  {\end{list}}%

\newenvironment{prf*}[1][Proof]{%
  \begin{proof}[\bf #1]
    \setcounter{equation}{0}
    } {\end{proof} }
\newcommand{\proofoftag}[2][:]{(#2)#1}
\newcommand{\proofofimp}[3][:]{\mbox{\eqclbl{#2}$\,\Rightarrow\,$\eqclbl{#3}#1}}

\newcommand{\pgref}[1]{\ref{#1}}
\newcommand{\thmref}[2][Theorem~]{#1\pgref{thm:#2}}
\newcommand{\corref}[2][Corollary~]{#1\pgref{cor:#2}}
\newcommand{\prpref}[2][Proposition~]{#1\pgref{prp:#2}}
\newcommand{\lemref}[2][Lemma~]{#1\pgref{lem:#2}}
\newcommand{\rmkref}[2][Remark~]{#1\pgref{rmk:#2}}
\newcommand{\secref}[2][Section~]{#1\ref{sec:#2}}
\renewcommand{\eqref}[1]{(\pgref{eq:#1})}
\newcommand{\thmcite}[2][?]{\cite[Thm.~#1]{#2}}
\newcommand{\corcite}[2][?]{\cite[Cor.~#1]{#2}}
\newcommand{\prpcite}[2][?]{\cite[Prop.~#1]{#2}}

\numberwithin{equation}{thm}

  \newcommand{\deq}{\:=\:} \newcommand{\dqis}{\:\qis\:}
\newcommand{\qis}{\simeq} \renewcommand{\le}{\leqslant}
 \newcommand{\Rhat}{\widehat{R}}

 \newcommand{\hev}[1]{\eta^{#1}}
\newcommand{\dptR}{\operatorname{depth}R}
\newcommand{\dpt}[2][R]{\operatorname{depth}_{#1}#2}
\newcommand{\id}[2][R]{\operatorname{id}_{#1}#2}
\newcommand{\pd}[2][R]{\operatorname{pd}_{#1}#2}
\newcommand{\Gid}[2][R]{\operatorname{Gid}_{#1}#2}
\newcommand{\Gpd}[2][R]{\operatorname{Gpd}_{#1}#2}
\newcommand{\RHom}[3][R]{\operatorname{\mathbf{R}Hom}_{#1}(#2,#3)}
\newcommand{\tp}[3][R]{\nobreak{#2\otimes_{#1}#3}}
\newcommand{\tpp}[3][R]{(\tp[#1]{#2}{#3})}
\newcommand{\Ltp}[3][R]{\nobreak{#2\otimes_{#1}^{\mathbf{L}}#3}}
\newcommand{\Ltpp}[3][R]{(\Ltp[#1]{#2}{#3})}
\newcommand{\catb}{\sqsubset\mspace{-13mu}\sqsupset}
\newcommand{\Cat}[2]{{\mathsf{#2}}(#1)}
\newcommand{\Catsup}[3]{{\mathsf{#2}}^{\text{\upshape #3}}(#1)}

\newcommand{\Catsupsub}[4]{{\mathsf{#2}}^{\text{\upshape
#3}}_{#4}(#1)} \newcommand{\D}[1][R]{\Cat{#1}{D}}
\newcommand{\Df}[1][R]{\Catsup{#1}{D}{f}}
\newcommand{\Dfb}[1][R]{\Catsupsub{#1}{D}{f}{\catb}}
\newcommand{\A}[1][R]{\Cat{#1}{A}} \newcommand{\B}[1][R]{\Cat{#1}{B}}
\newcommand{\Agross}[1][R]{\Cat{#1}{\hat{A}}}
\newcommand{\Bgross}[1][R]{\Cat{#1}{\hat{B}}}
\newcommand{\dcmca}[1]{\cite[#1]{dcmca}}
\newcommand{\aclassmap}[2]{\alpha_{#1}^{#2}}
\newcommand{\bclassmap}[2]{\beta_{#1}^{#2}}

\title[Homological dimensions of derived Hom
  complexes]{Homological dimensions\\ of derived Hom
  complexes}

\author[L.W.\ Christensen]{Lars Winther Christensen} %
\address{Texas Tech University, Lubbock, TX 79409, U.S.A.}
\email{lars.w.christensen@ttu.edu}
\urladdr{https://larswinther.github.io/homepage/}

\author[Soto Levins]{Andrew J. Soto Levins} %
\address{Texas Tech University, Lubbock, TX 79409, U.S.A.}
\email{ansotole@ttu.edu}
\urladdr{https://sites.google.com/view/andrewjsotolevins/home}

\begin{document}

\begin{abstract}
  Let $R$ be a commutative noetherian local ring and $M$ and $N$ be
  $R$-complexes with finitely generated homology. We prove that
  formulas for the homological dimensions of the derived Hom complex
  $\RHom{M}{N}$, in terms of the homological dimensions of $M$ and
  $N$, hold without the \emph{a priori} assumptions on the dimensions
  of $M$ or $N$ present in classic results.
\end{abstract}

\date{28 August 2026}

\keywords{Derived Hom complex, derived reflexive complex, perfect
  complex}

\subjclass{13D05; 13D07}

\maketitle

%%%%%%%%%%%%%%%%%%%%%%%%%%%%%%%%%%%%%%%%%%%%%%%%%%%%%%%%%%%%%%%%%%%%%%

\section*{Introduction}

\noindent
Throughout this paper, $R$ is a commutative noetherian local ring, and
$\Dfb$ denotes the full subcategory of the derived category, $\D$,
whose objects are complexes with bounded and degreewise finitely
generated homology. For complexes $M$ and $N$, in $\Dfb$ Foxby
\thmcite[4.1(a)]{HBF77b} proved the equality
\begin{equation}
  \label{eq:id}
  \id{\RHom{M}{N}} \deq \pd{M} + \id{N} \:;
\end{equation}
in particular, if two of the quantities $\id{\RHom{M}{N}}$, $\pd{M}$,
and $\id{N}$ are finite, then so is the third. Part (b) of the same
theorem asserts that the equality
\begin{equation}
  \label{eq:pd}
  \pd{\RHom{M}{N}} \deq \id{M} - \dpt{N}
\end{equation}
holds under the assumption that $\id{N}$ is finite. That is, under
that assumption, $\pd{\RHom{M}{N}}$ is finite if and only if $\id{M}$
is finite. In this paper we remove the \emph{a priori} assumption on
$N$ and prove, in \corref{perfect-I}, a genuine two-out-of-three
statement. That is to say, even though the quantity $\id{N}$ does not
show in equality \eqref{pd}, finiteness of $\pd{\RHom{M}{N}}$ and
$\id{M}$ forces it to be finite, and then the equality holds. Similar
improvements are made to other results from \cite{HBF77b,hha}; they
are also found in \secref{perfect}. Partial results in this
direction---imposing conditions on the homology of
$\RHom{M}{N}$---have recently been obtained by Dey, Holanda, and
Miranda-Neto \cite{DHM-24}, by Holanda, Jorge-P\'erez, and
Mendoza-Rubio \cite{HJM-25}, by Kimura \cite{Kimura:2026}, and by
Martins, Mendoza-Rubio, and Nason \cite{Martins/Mendoza/Nason:2026}.

Our proofs proceed via stability results for Gorenstein homological
dimensions. While it may not be explicitly stated in the literature,
one can piece together the equivalence of the following conditions for
a complex $M$ in $\Dfb$:
\begin{eqc}
\item $\Gpd{M}$ is finite.
\item $\RHom{M}{R}$ belongs to $\Dfb$ and $\Gpd{\RHom{M}{R}}$ is
  finite.
\end{eqc}
The key step towards our main results is to remove the boundedness
assumption on $\RHom{M}{R}$ in part \eqclbl{ii}, which is done in
\lemref{G-perfect}.

%%%%%%%%%%%%%%%%%%%%%%%%%%%%%%%%%%%%%%%%%%%%%%%%%%%%%%%%%%%%%%%%%%%%%%

\section{Stability of Auslander and Bass Categories revisited}

\noindent
Stability properties of Auslander and Bass Categories were studied by
Avramov and Foxby \cite{LLAHBF97} and have been revisited by several
authors since then.  We start by recalling some notation and
terminology from \dcmca{10.3} and, in general, we refer the reader to
\cite{dcmca} for any unexplained notation and terminology.

Assume that $R$ has a dualizing complex $D$ and consider the
adjunction,
\begin{equation*}
  \xymatrix@C=5.1pc{
    \D \ar@<2.5pt>[r]^-{\Ltp{D}{-}} & 
    \D\:. \ar@<2.5pt>[l]^-{\RHom{D}{-}}
  }
\end{equation*}  
For an $R$-complex $M$, the unit and counit of this adjunction are the
natural morphisms
\begin{equation*}
  \aclassmap{D}{M}\colon M\longrightarrow \RHom{D}{\Ltp{D}{M}}
\end{equation*}
and
\begin{equation*}
  \bclassmap{D}{M}\colon \Ltp{D}{\RHom{D}{M}}\longrightarrow M \:.
\end{equation*}
The gross Auslander Category of $R$ is the following full subcategory
of $\D$:
\begin{equation*}
  \Agross \deq \{M\in\D \mid \text{the unit morphism } \aclassmap{D}{M} \text{ is an isomorphism}\}
\end{equation*}
and the Auslander Category of $R$ is
\begin{equation*}
  \A \deq \{M\in\Agross \mid M \text{ and } \Ltp{D}{M} \text{ have bounded homology}\} \:.
\end{equation*}
The gross Bass Category of $R$ is the following full subcategory of
$\D$:
\begin{equation*}
  \Bgross \deq \{M\in\D\mid \text{the counit morphism } \bclassmap{D}{M} \text{ is an isomorphism}\}
\end{equation*}
and the Bass Category of $R$ is
\begin{equation*}
  \B \deq \{M\in\Bgross \mid M \text{ and } \RHom{D}{M} \text{ have bounded homology}\} \:.
\end{equation*}

\begin{ipg}
  Assume that $R$ has a dualizing complex $D$. For every complex $M$
  in $\Df$ the biduality morphism in $\D$,
  \begin{equation*}
    M \longrightarrow \RHom{\RHom{M}{D}}{D}
  \end{equation*}
  is an isomorphism; this is known as Grothendieck Duality. We recall
  from \dcmca{18.2.3, 19.1.7, 19.1.12, 19.1.13, 19.4.27} some
  properties used throughout the paper.
  \begin{itemlist}
  \item $M$ belongs to $\Dfb$ if and only if $\RHom{M}{D}$ belongs to
    $\Dfb$
  \end{itemlist}
  For a complex $M$ in $\Dfb$ one further has:
  \begin{itemlist}
  \item $\pd{M}$ is finite if and only if $\id{\RHom{M}{D}}$ is
    finite.
  \item $\Gpd{M}$ is finite if and only if $\Gid{\RHom{M}{D}}$ is
    finite, equivalently,\newline $M$ belongs to $\A$ if and only if
    $\RHom{M}{D}$ belongs to $\B$.
  \end{itemlist}
\end{ipg}

The next theorem is propelled by a result of Avramov, Iyengar, and
Lipman \thmcite[3.3]{AIL-10}; in the proof it is invoked through
\dcmca{19.4.19}.

\begin{thm}
  \label{thm:AB}
  Assume that $R$ has a dualizing complex.
  \begin{prt}
    
  \item For a complex $M$ in $\Dfb$ the following conditions are
    equivalent.
    \begin{eqc}
    \item $\Gpd{M}$ is finite.
    \item $M$ belongs to $\A$.
    \item $M$ belongs to $\Agross$.
    \item $\RHom{M}{D}$ belongs to $\Bgross$.
    \end{eqc}
  \item For a complex $N$ in $\Dfb$ the following conditions are
    equivalent.
    \begin{eqc}
    \item $\Gid{N}$ is finite.
    \item $N$ belongs to $\B$.
    \item $N$ belongs to $\Bgross$.
    \item $\RHom{N}{D}$ belongs to $\Agross$.
    \end{eqc}
  \end{prt}
\end{thm}

\begin{prf*}
  Let $D$ be a dualizing complex for $R$.
  
  \proofoftag{a} Conditions \eqclbl{i}--\eqclbl{iii} are equivalent by
  \dcmca{19.1.12, 19.4.19}. By \dcmca{19.1.7} it follows from
  \eqclbl{ii} that the complex $N = \RHom{M}{D}$ belongs to $\B$. To
  prove that \eqclbl{iv} implies \eqclbl{iii} it suffices by
  Grothendieck Duality to show that the complex $\RHom{N}{D}$ belongs
  to $\Agross$. There is a commutative diagram in $\D$,
  \begin{equation*}
    \xymatrix{
      \RHom{N}{D} \ar[r]^-{\RHom{\bclassmap{D}{N}}{D}} _-\qis
      \ar[d]^-{\aclassmap{D}{\RHom{N}{D}}}& \RHom{\Ltp{D}{\RHom{D}{N}}}{D}
      \ar[d]_-\qis\\
      \RHom{D}{\Ltp{D\!}{\!\RHom{N}{D}}} \ar[r]^-{\eta} & \mspace{-2mu} \RHom{D}{\RHom{\RHom{D}{N}}{D}}
    }
  \end{equation*}
  where $\eta = \RHom{D}{\hev{DND}}$ is an isomorphism by
  \dcmca{12.3.26(b)}. The diagram now shows that
  $\aclassmap{D}{\RHom{N}{D}}$ is an isomorphism, as desired.

  \proofoftag{b} Per \dcmca{19.1.13} the implications
  \eqclbl{i}$\,\Rightarrow\,$\eqclbl{ii}$\,\Rightarrow\,$\eqclbl{iii}
  are standard.

  \proofofimp{iii}{iv} Set $M = \RHom{N}{D}$. If $N \qis \RHom{M}{D}$
  belongs to $\Bgross$ then $M$ belongs to $\Agross$ by part (a).

  \proofofimp{iv}{i} The complex $M = \RHom{N}{D}$ belongs to $\Dfb$,
  so it has finite Gorenstein projective dimension by part (a), whence
  $N$ has finite Gorenstein injective dimension.
\end{prf*}

Recall that a complex $M$ in $\Dfb$ is called \emph{perfect} if it is
isomorphic in $\D$ to a bounded complex of finitely generated free
$R$-modules, equivalently $\pd{M}$ is finite. We recall a few facts
about perfect complexes that will be used throughout.

\begin{ipg}
  \label{P}
  Let $M$ be a perfect $R$-complex; one has
  \begin{equation}
    \label{eq:ABu}
    -\inf{\RHom{M}{R}} \deq \pd{M} \deq \dptR - \dpt{M} \:,
  \end{equation}
  and $\RHom{M}{R}$ is itself a perfect complex; see \dcmca{12.3.20,
    19.4.14}. The second equality above is the Auslander--Buchsbaum
  Formula from \cite{MAsDBc57}, see also \dcmca{16.4.2}.
\end{ipg}

The next result removes the boundedness assumption on $M$ in results
of Christensen and Holm \prpcite[4.1]{LWCHHl09} and Sather-Wagstaff
\prpcite[4.1]{SSW08}.

\enlargethispage*{\baselineskip}
\begin{prp}
  \label{prp:AB}
  Assume that $R$ has a dualizing complex and let $M$ and $N$ be
  complexes in $\Df$. If $N$ is perfect and not acyclic, then the next
  assertions hold.
  \begin{prt}
  \item The following conditions are equivalent.
    \begin{eqc}
    \item $M$ belongs to $\A$.
    \item $\Ltp{N}{M}$ belongs to $\A$.
    \item $\RHom{N}{M}$ belongs to $\A$.
    \end{eqc}
  \item The following conditions are equivalent.
    \begin{eqc}
    \item $M$ belongs to $\B$.
    \item $\Ltp{N}{M}$ belongs to $\B$.
    \item $\RHom{N}{M}$ belongs to $\B$.
    \end{eqc}
  \end{prt}
\end{prp}

\begin{prf*}
  First recall that if $M$ has bounded homology, then so does
  $\Ltp{N}{M}$; see e.g.\ \dcmca{15.4.3}. On the other hand, if
  $\Ltp{N}{M}$ has bounded homology, then $M$ has bounded homology by
  Iversen's amplitude inequality from \cite{BIv77} as improved by
  Foxby and Iyengar \thmcite[3.1]{HBFSIn03}, see also \dcmca{18.5.10}.

  Further, there are isomorphisms, see \dcmca{12.3.26(a)},
  \begin{align*}
    \RHom{N}{M} %%
    & \dqis \RHom{\RHom{\RHom{N}{R}}{R}}{M} \\
    & \dqis \Ltp{\RHom{N}{R}}{\RHom{R}{M}} \\
    & \dqis \Ltp{\RHom{N}{R}}{M} \:.
  \end{align*}
  As $\RHom{N}{R}$ is a perfect complex, $M$ has bounded homology if
  and only if $\RHom{N}{M}$ has bounded homology.
  
  The assertions now follow from \prpcite[4.1]{LWCHHl09} and
  \prpcite[4.1]{SSW08}.
\end{prf*}

\begin{cor}
  \label{cor:AB}
  Assume that $R$ has a dualizing complex and let $M$ and $N$ be
  complexes in $\Df$. If $N$ belongs to $\Dfb$, is not acyclic, and
  has finite injective dimension, then the next assertions hold.
  \begin{prt}
  \item $M$ belongs to $\A$ if and only if $\,\RHom{M}{N}$ belongs to
    $\B$.
  \item $M$ belongs to $\B$ if and only if $\,\RHom{M}{N}$ belongs to
    $\A$.
  \end{prt}
\end{cor}

\begin{prf*}
  Let $D$ be a dualizing complex for $R$; by Grothendieck Duality one
  has
  \begin{align*}
    \RHom{M}{N} & \dqis \RHom{M}{\RHom{\RHom{N}{D}}{D}} \\
                & \dqis{\RHom{\Ltp{\RHom{N}{D}}{M}}{D}} \:.
  \end{align*}
  The complex $\RHom{N}{D}$ is perfect, so the assertions are now
  immediate from \prpref{AB}.
\end{prf*}

%%%%%%%%%%%%%%%%%%%%%%%%%%%%%%%%%%%%%%%%%%%%%%%%%%%%%%%%%%%%%%%%%%%%%% 

\section{Derived reflexive complexes}
\label{sec:reflexive}

\noindent A complex $M$ in $\Dfb$ is called \emph{derived reflexive}
if $\RHom{M}{R}$ belongs to $\Dfb$ and the biduality morphism
$M \to \RHom{\RHom{M}{R}}{R}$ is invertible; see \dcmca{10.2.4}. Every
perfect $R$-complex is derived reflexive, see \dcmca{19.4.14}. We
recall a few facts about these complexes that will be used throughout
the paper.
  
\begin{ipg}
  \label{R}
  Let $M$ be a complex in $\Dfb$; it is derived reflexive if and only
  if it has finite Gorenstein projective dimension, see
  \dcmca{19.4.13}, and in that case one has
  \begin{equation}
    \label{eq:ABr}
    -\inf{\RHom{M}{R}} \deq \Gpd{M} \deq \dptR - \dpt{M} \:,
  \end{equation}
  and $\RHom{M}{R}$ is itself a derived reflexive complex, see
  \dcmca{19.4.14}.  The second equality above is the
  Auslander--Bridger Formula from \cite{MAsMBr69}, see also
  \dcmca{19.4.25}. Finally, if $R$ has a dualizing complex, then $M$
  is derived reflexive if and only if it belongs to the Auslander
  Category $\A$, see \dcmca{19.1.12}.
\end{ipg}

In several proofs the next remark allows us to assume the existence of
a dualizing $R$-complex.

\begin{rmk}
  \label{rmk:completion}
  For complexes $M$ and $N$ in $\Dfb$ the base changed complexes
  $\tp{\Rhat}{M}$ and $\tp{\Rhat}{N}$ belong to $\Dfb[\Rhat]$, see
  \dcmca{18.3.2}, and the complexes
  \begin{equation*}
    \tp{\Rhat}{\RHom{M}{N}} \dqis \RHom[\Rhat]{\tp{\Rhat}{M}}{\tp{\Rhat}{N}} 
  \end{equation*}
  and
  \begin{equation*}
    \tp{\Rhat}{\Ltpp{M}{N}} \dqis  \Ltp[\Rhat]{\tpp{\Rhat}{M}}{\tpp{\Rhat}{N}}
  \end{equation*}
  are in $\Df[\Rhat]$; see \dcmca{12.2.5, 12.2.10, 12.3.29,
    12.3.31(a)}.
    
  (a) For every complex $X$ in $\Df$ the projective and Gorenstein
  projective dimensions of $X$ over $R$ and $\tp{\Rhat}{X}$ over
  $\Rhat$ are simultaneously finite, see \dcmca{15.4.21, 17.4.28,
    19.3.15, 19.4.3}. In particular, $X$ is a derived reflexive or
  perfect $R$-complex if and only if $\tp{\Rhat}{X}$ is a derived
  reflexive or perfect $\Rhat$-complex.
  
  (b) For every complex $X$ in $\Df$ there is an inequality
  $\id[\Rhat]{\tpp{\Rhat}{X}} \le \id{X}$, see \dcmca{13.2.7,
    15.4.16(b), 17.3.20}; for $X$ in $\Dfb$ equality holds, see
  \dcmca{18.3.10}.
\end{rmk}

The apparent lack of duality between \eqref{pd} and the formula given
below for $\Gpd{\RHom{N}{M}}$ is easily clarified: While the dual
invariant to the infimum is the supremum, the infimum of a complex in
$\Df$ agrees per \dcmca{16.2.5(a)} with the width invariant which is,
indeed, dual to the depth.

\begin{prp}
  \label{prp:P-Gpd}
  Let $M$ and $N$ be complexes in $\Df$. If $N$ is perfect and not
  acyclic, then the following conditions are equivalent.
  \begin{eqc}
  \item $M$ is derived reflexive.
  \item $\Ltp{N}{M}$ is derived reflexive.
  \item $\RHom{N}{M}$ is derived reflexive.
  \end{eqc}
  When these conditions are satisfied, there are equalities,
  \begin{equation*}
    \Gpd{\Ltpp{N}{M}} \deq \pd{N} + \Gpd{M}
  \end{equation*}
  and
  \begin{equation*}
    \Gpd{\RHom{N}{M}} \deq \Gpd{M} - \inf{N} \:.
  \end{equation*}
\end{prp}

\begin{prf*}
  One can per \rmkref{completion} assume that $R$ is complete and, in
  particular, that it has a dualizing complex. The equivalence of
  conditions \eqclbl{i}--\eqclbl{iii} is now evident from
  \prpref{AB}(a).

  Assuming that conditions \eqclbl{i}--\eqclbl{iii} hold, one has
  \begin{align*}
    \Gpd{\Ltpp{N}{M}} %%
    & \deq \dptR -\dpt{\Ltpp{N}{M}} \\
    & \deq \dptR - (\dpt{N} + \dpt{M} - \dptR) \\   
    & \deq \pd{N} + \Gpd{M}   
  \end{align*}
  and
  \begin{align*}
    \Gpd{\RHom{N}{M}} %%
    & \deq \dptR -\dpt{\RHom{N}{M}} \\
    & \deq \dptR - (\dpt{M} + \inf{N}) \\
    & \deq \Gpd{M} - \inf{N} \:.   
  \end{align*}
  Here the first, fourth, and sixth equalities holds by the
  Auslander--Bridger Formula \eqref{ABr}. The second equality holds by
  \dcmca{16.3.1}, the third holds by the Auslander--Bridger and
  Auslander--Buchsbaum Formulas \eqref{ABr} and \eqref{ABu}. Finally,
  the fifth equality holds by \dcmca{16.2.25}.
\end{prf*}

\pagebreak
\begin{prp}
  \label{prp:P-Gid}
  Assume that $R$ has a dualizing complex and let $M$ and $N$ be
  complexes in $\Df$. If $N$ is perfect, then the following conditions
  are equivalent.
  \begin{eqc}
  \item $\Gid{M}$ is finite and $M$ belongs to $\Dfb$.
  \item $\Gid{\Ltpp{N}{M}}$ is finite and $\Ltp{N}{M}$ belongs to
    $\Dfb$.
  \item $\Gid{\RHom{N}{M}}$ is finite.
  \end{eqc}
  When these conditions are satisfied, there are equalities,
  \begin{equation*}
    \Gid{\Ltpp{N}{M}} \deq \Gid{M} - \inf{N}
  \end{equation*}
  and
  \begin{equation*}
    \Gid{\RHom{N}{M}} \deq \Gid{M} + \pd{N} \:.
  \end{equation*}
\end{prp}

\begin{prf*}
  The equivalence of the three conditions is a consequence of
  \dcmca{19.1.13} and \prpref{AB}(b).

  Assuming that conditions \eqclbl{i}--\eqclbl{iii} hold, one has
  \begin{align*}
    \Gid{\Ltpp{N}{M}} %%
    & \deq \dpt{R} - \inf{\Ltpp{N}{M}} \\
    & \deq \dptR - (\inf{N} + \inf{M}) \\          
    & \deq \Gid{M} - \inf{N}   
  \end{align*}
  and
  \begin{align*}
    \Gid{\RHom{N}{M}} %%
    & \deq \dptR -\inf{\RHom{N}{M}} \\
    & \deq \dptR - (\inf{M} + \dpt{N} - \dpt{R}) \\  
    & \deq \dptR - \inf{M} + (\dptR - \dpt{N}) \\
    & \deq \Gid{M} + \pd{N} \:.   
  \end{align*}
  The first, third, and fourth equalities hold by the Bass Formula for
  Gorenstein injective dimension, see \dcmca{19.2.40}, the second is
  standard, see \dcmca{16.2.10}. The fifth equality holds by
  \dcmca{16.3.5(b)}, and the last equality holds by the
  Auslander--Buchsbaum Formula \eqref{ABu} and another application of
  \dcmca{19.2.40}.
\end{prf*}

The equivalence of \eqclbl{i} and \eqclbl{ii} in the next lemma is
known from \thmcite[2.1]{AIL-10}. The difference between conditions
\eqclbl{ii} and \eqclbl{iii} is that the latter does not require the
homology of the complex $\RHom{M}{R}$ to be bounded below.

\begin{lem}
  \label{lem:G-perfect}
  Let $M$ be a complex in $\Dfb$; the next conditions are equivalent.
  \begin{eqc}
  \item $M$ is derived reflexive.
  \item $\RHom{M}{R}$ is derived reflexive.
  \item $\Gpd{\RHom{M}{R}}$ is finite.
  \end{eqc}
\end{lem}

\begin{prf*}
  Per \ref{R} the implications
  \eqclbl{i}$\,\Rightarrow\,$\eqclbl{ii}$\,\Rightarrow\,$\eqclbl{iii}
  are known.  To finish the proof, assume that $\Gpd{\RHom{M}{R}}$ is
  finite. Per \rmkref{completion} one can assume that $R$ has a
  dualizing complex; call it $D$. By \dcmca{19.1.12} the complex
  $\RHom{M}{R}$ belongs to $\Agross$. By \dcmca{19.1.2} the
  isomorphisms
  \begin{equation*}
    \RHom{M}{R} \dqis \RHom{M}{\RHom{D}{D}} \dqis \RHom{D}{\RHom{M}{D}}
  \end{equation*}
  now show that $\RHom{M}{D}$ belongs to $\Bgross$, whence $M$ is
  derived reflexive by \thmref{AB}(a).
\end{prf*}

\pagebreak
\begin{thm}
  \label{thm:perfect-G}
  Let $M$ and $N$ be complexes in $\Dfb$. If $N$ is perfect and not
  acyclic, then the following conditions are equivalent.
  \begin{eqc}
  \item $M$ is derived reflexive.
  \item $\RHom{M}{N}$ is derived reflexive.
  \item $\Gpd{\RHom{M}{N}}$ is finite.
  \end{eqc}
  When these conditions are satisfied, one has
  \begin{equation*}
    \Gpd{\RHom{M}{N}} \deq \pd{N} - \inf{M} \:.
  \end{equation*}
\end{thm}

\begin{prf*}
  The implication \eqclbl{i}\,$\,\Rightarrow\,$\,\eqclbl{ii} follows
  from \lemref{G-perfect} by induction on the number of non-zero
  modules in the complex of free $R$-modules isomorphic to $N$ in
  $\D$. Per \rmkref{completion} it is also a special case of
  \corcite[5.6]{LWCHHl09}. The implication
  \eqclbl{ii}$\,\Rightarrow\,$\eqclbl{iii} is standard. To see that
  \eqclbl{iii} implies \eqclbl{i}, assume $\Gpd{\RHom{M}{N}}$ is
  finite. Since the following isomorphisms hold,
  \begin{align*}
    \RHom{M}{N} & \dqis \RHom{M}{\RHom{\RHom{N}{R}}{R}} \\
                & \dqis \RHom{\Ltp{\RHom{N}{R}}{M}}{R} \:,
  \end{align*}
  and the complex $\RHom{N}{R}$ is perfect, the complex $M$ is derived
  reflexive by \lemref{G-perfect} and \prpref{P-Gpd}.

  Assuming that conditions \eqclbl{i}--\eqclbl{iii} hold, one has
  \begin{align*}
    \Gpd{\RHom{M}{N}} &\deq \dpt{R} - \dpt{\RHom{M}{N}} \\
                      &\deq \dpt{R} - (\dpt{N} + \inf{M}) \\
                      &\deq \pd{N} - \inf{M} \:.
  \end{align*}
  The first equality is the Auslander--Bridger Formula \eqref{ABr},
  the second holds by \dcmca{16.2.25}, and the third is the
  Auslander--Buchsbaum Formula \eqref{ABu}.
\end{prf*}

\begin{cor}
  \label{cor:perfect-GI}
  Assume that $R$ has a dualizing complex and let $M$ and $N$ be
  complexes in $\Dfb$. If $M$ is not acyclic and has finite injective
  dimension, then the following conditions are equivalent.
  \begin{eqc}
  \item $\Gid{N}$ is finite
  \item $\RHom{M}{N}$ is derived reflexive.
  \item $\Gpd{\RHom{M}{N}}$ is finite.
  \end{eqc}
  When these conditions are satisfied, one has
  \begin{equation*}
    \Gpd{\RHom{M}{N}} \deq \id{M} - \dpt{N} \:.
  \end{equation*}
\end{cor}

\begin{prf*}
  The implication \eqclbl{ii}$\,\Rightarrow\,$\eqclbl{iii} is
  standard. Now let $D$ be a dualizing complex for $R$ and note that
  by Grothendieck Duality, there are isomorphisms,
  \begin{align*}
    \RHom{M}{N} & \dqis \RHom{M}{\RHom{\RHom{N}{D}}{D}} \\
                & \dqis \RHom{\RHom{N}{D}}{\RHom{M}{D}} \:.
  \end{align*}
  The complex $\RHom{M}{D}$ is perfect.  Thus, if \eqclbl{i} holds,
  then $\RHom{N}{D}$ is derived reflexive, whence $\RHom{M}{N}$ is
  derived reflexive by \thmref{perfect-G}.  If \eqclbl{iii} holds,
  then the complex $\RHom{N}{D}$ is derived reflexive by
  \thmref{perfect-G}, and so $\Gid{N}$ is finite.

  Assuming that conditions \eqclbl{i}--\eqclbl{iii} hold one has
  \begin{align*}
    \Gpd{\RHom{M}{N}} &\deq \dpt{R} - \dpt{\RHom{M}{N}} \\
                      &\deq \dpt{R} - (\dpt{N} + \inf{M}) \\
                      &\deq \id{M} - \dpt{N}\:.
  \end{align*}
  Here the first equality is the Auslander--Bridger Formula
  \eqref{ABr}, the second equality holds by \dcmca{16.2.25} and the
  third by the Bass Formula, see \dcmca{16.4.11}.
\end{prf*}

%%%%%%%%%%%%%%%%%%%%%%%%%%%%%%%%%%%%%%%%%%%%%%%%%%%%%%%%%%%%%%%%%%%%%% 

\section{Perfect complexes}
\label{sec:perfect}

\noindent
The equivalence of conditions \eqclbl{i} and \eqclbl{ii} in the next
lemma is known from \thmcite[4.1]{AIL-10}. As in \lemref{G-perfect},
the difference between conditions \eqclbl{ii} and \eqclbl{iii} is that
the latter does not require the homology of $\RHom{M}{R}$ to be
bounded below.

\begin{lem}
  \label{lem:perfect}
  Let $M$ be a complex in $\Dfb$; the next conditions are equivalent.
  \begin{eqc}
  \item $M$ is perfect.
  \item $\RHom{M}{R}$ is perfect.
  \item $\pd{\RHom{M}{R}}$ is finite.
  \end{eqc}
\end{lem}

\begin{prf*}
  It is clear that the implications
  \eqclbl{i}$\,\Rightarrow\,$\eqclbl{ii}$\,\Rightarrow\,$\eqclbl{iii}
  hold. To finish the proof, assume that $\pd{\RHom{M}{R}}$ is
  finite. By \lemref{G-perfect} the complex $M$ is derived
  reflexive. In particular, the complex $\RHom{M}{R}$ belongs to
  $\Dfb$, so it is perfect and hence so is
  $M \qis \RHom{\RHom{M}{R}}{R}$.
\end{prf*}

The next result improves \dcmca{16.4.19}, which we know from Foxby's
notes \cite{hha}.

\begin{thm}
  \label{thm:perfect-P}
  Let $M$ and $N$ be non-acyclic complexes in $\Dfb$. If two of the
  conditions $(1)$--$(3)$ below are satisfied, then so is the third
  condition, and $\RHom{M}{N}$ is a perfect complex with
  \begin{equation*}
    \pd{\RHom{M}{N}} \deq \pd{N} - \inf{M} \:.      
  \end{equation*}
  \begin{rqm}
  \item $M$ is perfect.
  \item $N$ is perfect.
  \item $\pd{\RHom{M}{N}}$ is finite.
  \end{rqm}
\end{thm}

\begin{prf*}
  If $M$ is perfect, then $\RHom{M}{N}$ belongs to $\Dfb$ by
  \dcmca{15.4.3}, and \dcmca{16.4.19} yields
  $\pd{\RHom{M}{N}} = \pd{N} - \inf{M}$. Thus, it suffices to show
  that $(1)$ follows from $(2)$ and $(3)$. Since the following
  isomorphisms hold
  \begin{align*}
    \RHom{M}{N} & \dqis \RHom{M}{\RHom{\RHom{N}{R}}{R}} \\
                & \dqis \RHom{\Ltp{\RHom{N}{R}}{M}}{R}\:,
  \end{align*}
  the complex $\Ltp{\RHom{N}{R}}{M}$ is perfect by \lemref{perfect},
  whence $M$ is perfect by \thmcite[4.2(a)]{HBF77b}, which is also
  \dcmca{16.4.17}.
\end{prf*}

\begin{rmk}
  \label{rmk:pd}
  Notice from the proof above that if the complex $N$ is derived
  reflexive and $\pd{\RHom{M}{N}}$ is finite, then $M$, $N$, and
  $\RHom{M}{N}$ are perfect. This was proved in
  \thmcite[1.1(1)]{Martins/Mendoza/Nason:2026} under the \emph{a
    priori} assumption that the complex $\RHom{M}{N}$ has bounded
  homology.
\end{rmk}

The next result improves \thmcite[4.1(b)]{HBF77b} which is also
\dcmca{16.4.34}.

\begin{cor}
  \label{cor:perfect-I}
  Let $M$ and $N$ be non-acyclic complexes in $\Dfb$. If two of the
  conditions $(1)$--$(3)$ below are satisfied, then so is the third
  condition, and $\RHom{M}{N}$ is a perfect complex with
  \begin{equation*}
    \pd{\RHom{M}{N}} \deq \id{M} - \dpt{N} \:.
  \end{equation*}
  \begin{rqm}
  \item $\id{M}$ is finite.
  \item $\id{N}$ is finite.
  \item $\pd{\RHom{M}{N}}$ is finite.
  \end{rqm}
\end{cor}

\begin{prf*}
  If $\id{N}$ is finite, then $\RHom{M}{N}$ belongs to $\Dfb$ by
  \dcmca{15.4.9}, and \dcmca{16.4.34} yields
  $\pd{\RHom{M}{N}} = \id{M} - \dpt{N}$. Thus, it suffices to show
  that $(2)$ follows from $(1)$ and $(3)$. Per \rmkref{completion} one
  can assume that $R$ has a dualizing complex, call it $D$. The
  following isomorphisms hold
  \begin{align*}
    \RHom{M}{N} & \dqis \RHom{M}{\RHom{\RHom{N}{D}}{D}} \\
                & \dqis \RHom{\RHom{N}{D}}{\RHom{M}{D}} \:.
  \end{align*}
  The complex $\RHom{M}{D}$ is perfect, so $\RHom{N}{D}$ is perfect by
  \thmref{perfect-P}, whence $\id{N}$ is finite.
\end{prf*}

\begin{rmk}
  \label{rmk:pd_id}
  Notice from the proof above and \rmkref{pd} that if $M$ has finite
  Gorenstein injective dimension and $\pd{\RHom{M}{N}}$ is finite,
  then $M$ and $N$ have finite injective dimension, and $\RHom{M}{N}$
  is perfect. This was proved in
  \thmcite[1.1(2)]{Martins/Mendoza/Nason:2026} under the \emph{a
    priori} assumption that $\RHom{M}{N}$ is in $\Dfb$.
\end{rmk}

For complexes $M$ and $N$ in $\Dfb$ it is shown in
\thmcite[4.2(a)]{HBF77b} that
\begin{equation*}
  \pd{\Ltpp{M}{N}} \deq \pd{M} + \pd{N}
\end{equation*}
holds without further assumptions, while
$\id{\Ltpp{M}{N}} = \id{M} - \inf{N}$ by part (b) of the same theorem
holds under the assumption that $N$ is perfect. The next result
improves \thmcite[4.2(b)]{HBF77b}, which is also \dcmca{16.4.35}.
\enlargethispage*{\baselineskip}

\begin{thm}
  \label{thm:perfect-I-tp}
  Let $M$ and $N$ be non-acyclic complexes in $\Dfb$. If two of the
  conditions \mbox{$(1)$--$(3)$} below are satisfied, then so is the
  third condition, and $\Ltp{M}{N}$ is a complex in $\Dfb$ of finite
  injective dimension:
  \begin{equation*}
    \id{\Ltpp{M}{N}} \deq \id{M} - \inf{N} \:.
  \end{equation*}
  \begin{rqm}
  \item $\id{M}$ is finite.
  \item $N$ is perfect.
  \item $\id{\Ltpp{M}{N}}$ is finite.
  \end{rqm}
\end{thm}

\begin{prf*}
  If $N$ is perfect, then $\Ltp{M}{N}$ belongs to $\Dfb$, see
  \dcmca{15.4.3}, and \dcmca{16.4.35} yields
  $\id{\Ltpp{M}{N}} = \id{M} - \inf{N}$. Thus, it suffices to show
  that $(2)$ follows from $(1)$ and $(3)$. Set
  $\widehat{M} = \tp{\Rhat}{M}$ and $\widehat{N} = \tp{\Rhat}{N}$; by
  \rmkref{completion} the $\Rhat$-complexes $\widehat{M}$ and
  $\Ltp[\Rhat]{\widehat{M}}{\widehat{N}}$ have finite injective
  dimension, and it suffices to show that $\widehat{N}$ is a perfect
  $\Rhat$-complex.
  
  Let $D$ be a dualizing complex for $\Rhat$. First we argue that the
  complex $\Ltp[\Rhat]{\widehat{M}}{\widehat{N}}$ has bounded
  homology.  Since $\widehat{M}$ belongs to the Bass category, see
  \dcmca{10.3.5}, there are isomorphisms,
  \begin{align*}
    \Ltp[\Rhat]{\widehat{M}}{\widehat{N}} %%
    & \dqis \Ltp[\Rhat]{\Ltpp[\Rhat]{D}{\RHom[\Rhat]{D}{\widehat{M}}}}{\widehat{N}} \\
    & \dqis \Ltp[\Rhat]{D}{\Ltpp[\Rhat]{\RHom[\Rhat]{D}{\widehat{M}}}{\widehat{N}}} \:.
  \end{align*}
  As $\RHom[\Rhat]{D}{\widehat{M}}$ is perfect, the complex
  $\Ltp[\Rhat]{\RHom[\Rhat]{D}{\widehat{M}}}{\widehat{N}}$ belongs to
  $\Dfb$, see \dcmca{15.4.9, 15.4.3}.  Since
  $\Ltp[\Rhat]{\widehat{M}}{\widehat{N}}$ belongs to $\Bgross$, again
  by \dcmca{10.3.5}, the complex
  $\Ltp[\Rhat]{\RHom[\Rhat]{D}{\widehat{M}}}{\widehat{N}}$ belongs to
  $\Agross$, see \dcmca{19.1.2}, and therefore to $\A$ by
  \thmref{AB}(a). Now another application of \dcmca{19.1.2} shows that
  $\Ltp[\Rhat]{\widehat{M}}{\widehat{N}}$ belongs to $\B$; in
  particular, it has bounded homology.  Now the dual complex
  \begin{equation*}
    \RHom[\Rhat]{\Ltp[\Rhat]{\widehat{M}}{\widehat{N}}}{D} \dqis \RHom[\Rhat]{\widehat{N}}{\RHom[\Rhat]{\widehat{M}}{D}}
  \end{equation*}
  is perfect, and since $\RHom[\Rhat]{\widehat{M}}{D}$ is perfect, so
  is $\widehat{N}$ by \thmref{perfect-P}.
\end{prf*}
\enlargethispage*{2\baselineskip}

\providecommand{\MR}{\relax\ifhmode\unskip\space\fi MR }
\providecommand{\MRhref}[2]{%
  \href{http://www.ams.org/mathscinet-getitem?mr=#1}{#2} }

\end{document}